\documentclass[12pt,reqno]{article}

\usepackage[usenames]{color}
\usepackage{amssymb}
\usepackage{graphicx}
\usepackage{amscd}

\usepackage{amsthm}
\newtheorem{theorem}{Theorem}

\newtheorem{proposition}[theorem]{Proposition}
\newtheorem{corollary}[theorem]{Corollary}

\theoremstyle{definition}

\newtheorem{example}[theorem]{Example}

\usepackage[colorlinks=true,
linkcolor=webgreen, filecolor=webbrown,
citecolor=webgreen]{hyperref}

\definecolor{webgreen}{rgb}{0,.5,0}
\definecolor{webbrown}{rgb}{.6,0,0}

\usepackage{color}

\usepackage{float}

\usepackage{graphics,amsmath,amssymb}
\usepackage{amsfonts}
\usepackage{latexsym}
\usepackage{epsf}

\newcommand{\seqnum}[1]{\href{http://www.research.att.com/cgi-bin/access.cgi/as/~njas/sequences/eisA.cgi?Anum=#1}{\underline{#1}}}

\begin{document}

\begin{center}
\vskip 1cm{\LARGE\bf An orthogonal polynomial coefficient formula for the Hankel transform} \vskip 1cm \large
Paul Barry\\
School of Science\\
Waterford Institute of Technology\\
Ireland\\
\href{mailto:pbarry@wit.ie}{\tt pbarry@wit.ie}
\end{center}
\vskip .2 in

\begin{abstract} We give an explicit formula for the Hankel transform of a regular sequence in terms of the coefficients of the associated orthogonal polynomials and the sequence itself. We apply this formula to some sequences of combinatorial interest, deriving interesting combinatorial identities by this means. Further insight is also gained into the structure of the Hankel transform of a sequence.
\end{abstract}

Let $\mu_n$ denote a sequence of numbers, with $\mu_0=1$. The Hankel transform \cite{Layman} of this sequence is the sequence $h_n$ of determinants
$$h_n=| \mu_{i+j} |_{0 \le i,j \le n}.$$
\noindent We shall assume that the sequence $\mu_n$ is \emph{regular} (or ``Catalan-like''), that is, $h_n \ne 0$ for any $n$.
Much research \cite{BRP, CRI, Ege, ExpPoly, Junod, Petkovic, RadEff, RadHermite, RadDer} has focused on finding closed form expressions for the Hankel transforms of important sequences. Some well-known sequences have Hankel transforms with simple formulas (such as our examples below), but in general Hankel transforms can have quite complex forms \cite{Ege}, or be totally resistant to current techniques of elucidation. The sequences we shall use in this note are susceptible to the following classical approach.
If \cite{Kratt, Kratt1} $\mu_n$ has a generating function  $g(x)$ expressible in the form
$$g(x)=\cfrac{\mu_0}{1-\alpha_0 x-
\cfrac{\lambda_1 x^2}{1-\alpha_1 x-
\cfrac{\lambda_2 x^2}{1-\alpha_2 x-
\cfrac{\lambda_3 x^2}{1-\alpha_3 x-\cdots}}}}$$ then
we have \cite{Kratt}
\begin{equation}\label{Kratt} h_n = \mu_0^{n+1} \lambda_1^n\lambda_2^{n-1}\cdots \lambda_{n-1}^2\lambda_n=\mu_0^{n+1}\prod_{k=1}^n
\lambda_k^{n+1-k}.\end{equation}
We can associate a sequence of monic orthogonal polynomials \cite{Chihara, Gautschi, Szego} $P_n(x)$ with the sequence $\mu_n$. The $\mu_n$ then correspond to the moments of this family, for the linear functional $\mathcal{L}$ defined by
$$\mathcal{L}(x^n)=\mu_n.$$ The family $\{P_n(x)\}$ satisfies a three term recurrence relation
$$P_{n+1}(x)=(x-\alpha_n) P_n(x)-\lambda_n P_{n-1}(x), \quad n \ge 0,$$ where $\lambda_0P_{-1}(x)=0$ and $P_1(x)=1$. 
We can construct the $P_n(x)$ explicitly as follows \cite{Ismail}. Let 
$$\Delta_{i,n}=\begin{vmatrix} \mu_0 & \mu_1 & \cdots & \mu_i \\
                             \mu_1 & \mu_2 & \cdots & \mu_{i+1}\\
                             \vdots & \vdots & \vdots & \vdots \\
                             \mu_{i-1} & \mu_i & \cdots & \mu_{2i-1}\\
                             \mu_n & \mu_{n+1} & \cdots & \mu_{n+1} \end{vmatrix}, \quad 
                             D_n(x)= 
                             \begin{vmatrix} \mu_0 & \mu_1 & \cdots & \mu_i \\
                             \mu_1 & \mu_2 & \cdots & \mu_{i+1}\\
                             \vdots & \vdots & \vdots & \vdots \\
                             \mu_{i-1} & \mu_i & \cdots & \mu_{2i-1}\\
                             1 & x & \cdots & x^n \end{vmatrix},$$
and let 
$$D_n= \Delta_{n,n}.$$ 
\noindent Then $$P_n(x)=\frac{D_n(x)}{D_{n-1}}$$ is the monic orthogonal sequence for $\mathcal{L}$. 
The orthogonality of the polynomials is equivalent to
$$\mathcal{L}(P_m P_n)=\lambda_1 \lambda_2 \cdots \lambda_n \delta_{m,n}.$$
In particular, $$\mathcal{L}(P_n^2)=\lambda_1 \lambda_2 \cdots \lambda_n.$$
\noindent However, we know that \cite{RadEff}
$$h_n = \prod_{k=0}^n \lambda_1 \lambda_2 \cdots \lambda_k.$$
Thus we have that
$$h_n=\prod_{k=0}^n \mathcal{L}(P_k^2).$$
We wish now to evaluate $\mathcal{L}(P_k^2)$. To do this, we let $(a_{n,k})$ be the matrix of coefficients of the polynomials $\{P_n(x)\}$. Since the degree of $P_n(x)$ is $n$, this is a lower triangular matrix, with $1$'s on the diagonal since the polynomials are monic. It is thus invertible (and if its elements are integers, so are those of the inverse), and we know that \cite{Riordan_OP} the first column of its inverse is composed of the $\mu_n$.
Using the Cauchy formula for the coefficients of the product of two polynomials, along with the fact that $\mathcal{L}$ operates on
polynomials as follows:
$$\mathcal{L}(\sum_{j=0}^m a_j x^j)=\sum_{j=0}^m a_j \mu_j,$$ we obtain the following formula:
\begin{equation}\label{EqL}\mathcal{L}(P_k^2)=\sum_{i=0}^{2k} \left(\sum_{j=0}^i a_{k,j} a_{k,i-j}\right)\mu_i.\end{equation}
\noindent Thus we obtain the following formula for the Hankel transform $h_n$ in terms of the elements of the coefficient array for the associated orthogonal polynomials:

\begin{equation}\label{EqH}h_n = \prod_{k=0}^n \sum_{i=0}^{2k} (\sum_{j=0}^i a_{k,j} a_{k,i-j})\mu_i.\end{equation}
We state this as a proposition. 
\begin{proposition} Let $\mu_n$ be a regular sequence, whose associated family of orthogonal polynomials $\{P_n(x)\}_{n \ge 0}$ has coefficient array $\left(a_{n,k}\right)_{n,k \ge 0}$. Then the Hankel transform $h_n$ of $\mu_n$ is given by 
$$h_n=\prod_{k=0}^n \sum_{i=0}^{2k} (\sum_{j=0}^i a_{k,j} a_{k,i-j})\mu_i.$$
\end{proposition}

\noindent Note also that the elements of this product, namely $\lambda_1 \lambda_2 \cdots \lambda_k$, are precisely the elements of
of $D$ where \cite{Triple, P_W, Woan} the Hankel matrix $H$ of $\mu_n$ has the $LDU$ decomposition
$$H_{\mu} = \left(\mu_{i+j}\right) = L D L^T,$$ where $L$ is the inverse of $\left(a_{i,j}\right)$. Thus the elements of $D$ are given by $$ \sum_{i=0}^{2n} (\sum_{j=0}^i a_{n,j} a_{n,i-j})\mu_i, \quad n=0,1,2,\ldots.$$ \noindent We state this as a corollary to the proposition.
\begin{corollary} Let $\mu_n$ be a regular sequence, with Hankel matrix $H=\left(\mu_{i+j}\right)_{i,j\ge 0}$. Then 
$$H = LDL^T,$$ where the elements $d_n$ of the diagonal matrix $D$ are expressible as 
$$d_k = \sum_{i=0}^{2k} (\sum_{j=0}^i a_{k,j} a_{k,i-j})\mu_i.$$
\end{corollary}
\section{Examples}
In this section, we give some examples, concentrating on ``simple'' sequences of combinatorial significance whose Hankel transforms are well known. In each case the coefficient array of the associated family of monic orthogonal polynomials is a Riordan array \cite{DeutschShap, SGWW, Spru}. These examples further explore links \cite{Riordan_OP, Meixner, Toda} between Riordan arrays, orthogonal polynomials and Hankel transforms. For those not familiar with Riordan arrays, we recall that 
the \emph{Riordan group} \cite{SGWW, Spru}, is a set of
infinite lower-triangular integer matrices, where each matrix is
defined by a pair of generating functions
$g(x)=1+g_1x+g_2x^2+\cdots$ and $f(x)=f_1x+f_2x^2+\cdots$ where
$f_1\ne 0$ \cite{Spru}. We assume in addition that $f_1=1$ in what follows. The associated matrix is the matrix whose
$i$-th column is generated by $g(x)f(x)^i$ (the first column being
indexed by 0). The matrix corresponding to the pair $g, f$ is
denoted by $(g, f)$ or $\cal{R}$$(g,f)$. The group law is then given
by
\begin{displaymath} (g, f)\cdot(h, l)=(g, f)(h, l)=(g(h\circ f), l\circ
f).\end{displaymath} The identity for this law is $I=(1,x)$ and the
inverse of $(g, f)$ is $(g, f)^{-1}=(1/(g\circ \bar{f}), \bar{f})$
where $\bar{f}$ is the compositional inverse of $f$ (defined by $f(\bar{f}(x))=\bar{f}(f(x))=x)$).
\begin{example} The binomial matrix (Pascal's triangle) which is the matrix with general term $\binom{n}{k}$ is the Riordan array 
$$\mathbf{B}=\left(\frac{1}{1-x},\frac{x}{1-x}\right).$$ It has inverse $\mathbf{B}^{-1}=\left((-1)^{n-k}\binom{n}{k}\right)$ given by $$\mathbf{B}^{-1}=\left(\frac{1}{1-x},\frac{x}{1-x}\right).$$ 
\end{example}
\noindent
 The \emph{exponential Riordan group} \cite
{DeutschShap, ProdMat}, is a set of
infinite lower-triangular integer matrices, where each matrix
is defined by a pair
of generating functions $g(x)=g_0+g_1x+g_2x^2+\cdots$ and
$f(x)=f_1x+f_2x^2+\cdots$ where $g_0 \ne 0$ and $f_1\ne 0$. In what follows, we shall assume
$$g_0=f_1=1.$$
The associated
matrix is the matrix
whose $i$-th column has exponential generating function
$g(x)f(x)^i/i!$ (the first column being indexed by 0). The
matrix corresponding to
the pair $f, g$ is denoted by $[g, f]$.  The group law is given by \begin{displaymath}
[g,
f]\cdot [h,
l]=[g(h\circ f), l\circ f].\end{displaymath} The identity for
this law is $I=[1,x]$ and the inverse of $[g, f]$ is $[g,
f]^{-1}=[1/(g\circ
\bar{f}), \bar{f}]$ where $\bar{f}$ is the compositional
inverse of $f$. 
\begin{example} The binomial matrix $\mathbf{B}=\left(\binom{n}{k}\right)$ is also an element of the 
exponential Riordan group. We have 
$$\mathbf{B}=\left[e^x, x\right] \quad \text{and} \quad \mathbf{B}^{-1}=\left[e^{-x},x\right].$$ 
\end{example}

We are now is a position to apply our proposition to some sequences of combinatorial interest. 
\begin{example} The Catalan numbers \seqnum{A000108}. 
It is well known that for the Catalan numbers $C_n=\frac{1}{n+1}\binom{2n}{n}$ we have $h_n=1$ for all $n$. This can be seen since  the generating function $g(x)=\frac{1-\sqrt{1-4x}}{2x}$ of the Catalan numbers can be expressed as 
$$g(x)=\cfrac{1}{1-x-\cfrac{x^2}{1-2x-\cfrac{x^2}{1-2x-\cfrac{x^2}{1-\cdots}}}}.$$ 
\noindent The corresponding monic orthogonal polynomials
$$P_n(x)=\sum_{k=0}^n (-1)^{n-k}\binom{n+k}{2k}x^k$$ have coefficient array $\left(\frac{1}{1+x},\frac{x}{(1+x)^2}\right)$. We have
$$P_{n+1}(x)=(x-2)P_n(x)-P_{n-1}(x),\quad P_0(x)=1,\quad P_1(x)=x-1.$$
Equation (\ref{EqH}) takes the form
$$\prod_{k=0}^n \sum_{i=0}^{2k} \sum_{j=0}^i (-1)^i \binom{k+j}{2j}\binom{k+i-j}{2(i-j)} C_i = 1.$$
In fact in this case, since $\lambda_i=1$, we obtain the identity
$$\sum_{i=0}^{2k} \sum_{j=0}^i (-1)^i \binom{k+j}{2j}\binom{k+i-j}{2(i-j)} C_i = 1.$$
\end{example}
\begin{example} The central binomial numbers \seqnum{A000984}.  The Hankel transform of the central binomial numbers $\binom{2n}{n}$ is given by
$h_n=2^n$. This can be seen since  the generating function $g(x)=\frac{1}{\sqrt{1-4x}}$ of the central binomial numbers can be expressed as
$$g(x)=\cfrac{1}{1-2x-\cfrac{2x^2}{1-2x-\cfrac{x^2}{1-2x-\cfrac{x^2}{1-\cdots}}}}.$$ 

\noindent The corresponding orthogonal polynomials $P_n(x)$ have coefficient array $\left(\frac{1-x}{1+x},\frac{x}{(1+x)^2}\right)$. In this case, we have
$$P_{n+1}(x)=(x-2)P_n(x)-P_{n-1}(x),\quad P_0(x)=1,\quad P_1(x)=x-2.$$ The general term $a_{n,k}$ of the Riordan array $\left(\frac{1-x}{1+x},\frac{x}{(1+x)^2}\right)$ is given by
$$a_{n,k}=\left(\binom{n+k}{n-k}+\binom{n+k-1}{n-k-1}\right)(-1)^{n-k}=\binom{n+k}{2k}\frac{2n+0^{n+k}}{n+k+0^{n+k}}(-1)^{n-k}.$$ Here we have used the notation $0^n$ to denote the sequence with term $1,0,0,0,\ldots$ and generating function $1$. Thus $0^n=\delta_{0,n}=[n=0]$.
With this value for $a_{n,k}$, we then obtain
$$\prod_{k=0}^n \sum_{i=0}^{2k} (\sum_{j=0}^i a_{k,j} a_{k,i-j})\binom{2i}{i} = 2^n.$$
Now $\lambda_1=2$, $\lambda_i=1$ for $i>1$, and so we also have
$$\sum_{i=0}^{2k} (\sum_{j=0}^i a_{k,j} a_{k,i-j})\binom{2i}{i}=2-0^k.$$ Thus we have, for instance,
$$\sum_{i=0}^{2k} \sum_{j=0}^i (-1)^i \binom{k+j}{2j}\binom{k+i-j}{2(i-j)}\frac{2k+0^{k+j}}{k+j+0^{k+j}}\frac{2k+0^{k+i-j}}{k+j+0^{k+i-j}}\binom{2i}{i}=2-0^k.$$
\end{example}
\begin{example} The large Schr\"oder numbers \seqnum{A006318}.  The large Schr\"oder numbers $S_n=\sum_{k=0}^n \binom{n+k}{2k}C_k$ have Hankel transform $h_n=2^{\binom{n+1}{2}}.$ This can be seen since  the generating function $g(x)$  of the large Schr\"oder numbers can be expressed as
$$g(x)=\cfrac{1}{1-2x-\cfrac{2x^2}{1-3x-\cfrac{2x^2}{1-3x-\cfrac{2x^2}{1-\cdots}}}}.$$ The corresponding monic orthogonal polynomials $P_n(x)$ satisfy
$$P_{n+1}(x)=(x-3)P_n(x)-2P_{n-1}(x),\quad P_0(x)=1,\quad P_1(x)=x-2,$$ with coefficient array equal to the Riordan array
$$\left(\frac{1}{1+2x},\frac{x}{1+3x+2x^2}\right).$$ For this array, we have
$$a_{n,k}=(-1)^{n-k}\sum_{j=0}^n \binom{n-j}{k}\binom{n+k}{j}.$$ With this value for $a_{n,k}$, we then have
$$h_n = 2^{\binom{n+1}{2}}=\prod_{k=0}^n \sum_{i=0}^{2k} (\sum_{j=0}^i a_{k,j} a_{k,i-j})S_i,$$ and
$$ \sum_{i=0}^{2k} (\sum_{j=0}^i a_{k,j} a_{k,i-j})S_i = 2^k.$$
Explicitly, this last equation gives us
$$\sum_{i=0}^{2k} \sum_{j=0}^i\sum_{l=0}^k \binom{k-l}{j}\binom{k+j}{l}\sum_{m=0}^k \binom{k-m}{i-j}\binom{k+i-j}{m} S_i=2^k.$$
\end{example}
\begin{example} The factorial numbers \seqnum{A000142}. The factorial numbers $n!$ have Hankel transform $h_n=\prod_{i=0}^n i!^2$. This can be seen since  the generating function $g(x)$  of the factorial numbers can be expressed as
$$g(x)=\cfrac{1}{1-x-\cfrac{x^2}{1-3x-\cfrac{4x^2}{1-5x-\cfrac{9x^2}{1-\cdots}}}}.$$

\noindent The corresponding monic orthogonal polynomials $P_n(x)$ satisfy
$$P_{n+1}(x)=(x-(2n+1))P_n(x)-n^2 P_{n-1}(x), \quad P_0(x)=1,\quad P_1(x)=x-1,$$ with coefficient array equal to the exponential Riordan array $$\left[\frac{1}{1+x},\frac{x}{1+x}\right]$$ with general term $a_{n,k}$ given by
$$a_{n,k}=\binom{n}{k}\frac{n!}{k!}(-1)^{n-k}.$$
With this value for $a_{n,k}$, we then obtain
$$\prod_{k=0}^n \sum_{i=0}^{2k} (\sum_{j=0}^i a_{k,j} a_{k,i-j})i! = \prod_{i=0}^n i!^2,$$ and
$$\sum_{i=0}^{2k} (\sum_{j=0}^i a_{k,j} a_{k,i-j})i!=k!^2.$$ Writing these expressions explicitly, we get
$$\prod_{k=0}^n \sum_{i=0}^{2k} (\sum_{j=0}^i (-1)^i \binom{k}{j}\binom{k}{i-j}\frac{k!^2}{j!(i-j)!})i! = \prod_{i=0}^n i!^2,$$
and
$$\sum_{i=0}^{2k} (\sum_{j=0}^i (-1)^i \binom{k}{j}\binom{k}{i-j}\frac{k!^2}{j!(i-j)!})i!=k!^2.$$ We deduce that
$$\sum_{i=0}^{2k} (\sum_{j=0}^i (-1)^i \binom{k}{j}\binom{k}{i-j}\frac{1}{j!(i-j)!})i!=1,$$ or
$$\sum_{i=0}^{2k} (\sum_{j=0}^i (-1)^i \binom{k}{j}\binom{k}{i-j}\binom{i}{j})=1.$$ Expanding the inner term in $i$ and $k$ as an array, we see that this expression represents the row sums of the array that begins
\begin{equation*}\begin{array}{ccccccccccc}1&&&&&&&&&&\\1&-2&2&&&&&&&&\\1&-4&10&-12&6&&&&&&\\1&-6&24&-56&78&-60&20&&&&\\1&-8&44&-152&346&-520&500&-280&70&& \\
1&-10&70&-320&1010&-2252&3560&-3920&2870&-1260&252\end{array}\end{equation*}

\end{example}
\begin{example} The derangement (or rencontres) numbers \seqnum{A000166}. 
It is well-known \cite{Layman} that a sequence and its binomial transform, or its inverse binomial transform, have the same Hankel transform. Thus the derangement numbers \cite{RadDer} 
$$d_n=\sum_{k=0}^n (-1)^{n-k}\binom{n}{k}k!,$$ which are the inverse binomial transform of the factorial numbers $n$!, have Hankel transform $h_n=\prod_{i=0}^n i!^2$. This may be seen from the continued fraction expansion of their generating function, given by 
$$\cfrac{1}{1-\cfrac{x^2}{1-2x-\cfrac{4x^2}{1-4x-\cfrac{9x^2}{1-\cdots}}}}. $$ In this case, the coefficient array of the associated orthogonal polynomials is equal to that of the factorial numbers, multiplied on the left by the  binomial matrix. We then obtain 
$$[e^x,x]\cdot \left[\frac{1}{1+x},\frac{x}{1+x}\right]=\left[\frac{e^x}{1+x},\frac{x}{1+x}\right],$$ as the Riordan array expression for the coefficient array of the associated orthogonal polynomials, with general term
$a_{n,k}$ given by 
$$a_{n,k}=\sum_{j=0}^n \binom{n}{j}\binom{j}{k}\frac{j!}{k!}(-1)^{j-k}.$$ 
\noindent With this value for $a_{n,k}$, we then deduce, for instance, that 
$$\sum_{i=0}^{2k} (\sum_{j=0}^i a_{k,j} a_{k,i-j})\sum_{l=0}^i \binom{i}{l}(-1)^{i-l}l!=k!^2.$$
\end{example}

\bigskip
\hrule
\bigskip
\noindent 2010 {\it Mathematics Subject Classification}:
Primary 42C05; Secondary
11B83, 11C20, 15B05, 15B36, 33C45.

\noindent \emph{Keywords:} Integer sequence, orthogonal polynomials, moments,
Riordan array, Hankel determinant, Hankel transform.

\bigskip
\hrule
\bigskip
\noindent Concerns sequences
\seqnum{A000108},
\seqnum{A000142},
\seqnum{A000166},
\seqnum{A000984},
\seqnum{A006318}.

\end{document}